\RequirePackage{ifpdf}
\ifpdf 
\documentclass[pdftex]{sigma}
\else
\documentclass{sigma}
\fi

\begin{document}

\allowdisplaybreaks

\renewcommand{\thefootnote}{$\star$}

\renewcommand{\PaperNumber}{005}

\FirstPageHeading

\ShortArticleName{Curvature Operators in Lorentzian Manifolds with
Large Isometry Groups}

\ArticleName{Algebraic Properties of Curvature Operators \\ in Lorentzian Manifolds with Large Isometry Groups}

\Author{Giovanni CALVARUSO~$^\dag$ and Eduardo GARC\'IA-R\'IO~$^\ddag$}

\AuthorNameForHeading{G.~Calvaruso and E.~Garc\'{\i}a-R\'{\i}o}

\Address{$^\dag$~Dipartimento di Matematica ``E.~De Giorgi'', Universit\`{a} del Salento, Lecce, Italy}
\EmailD{\href{mailto:giovanni.calvaruso@unisalento.it}{giovanni.calvaruso@unisalento.it}}

\Address{$^\ddag$~Faculty of Mathematics, University of  Santiago de Compostela,\\
 \hphantom{$^\ddag$}~15782 Santiago de Compostela, Spain}
\EmailD{\href{mailto:eduardo.garcia.rio@usc.es}{eduardo.garcia.rio@usc.es}}

\ArticleDates{Received October 01, 2009, in f\/inal form January 07, 2010;  Published online January 12, 2010}

\Abstract{Together with spaces of constant sectional curvature and products of a real line with a manifold of constant curvature, the socalled {\em Egorov spaces} and {\em $\varepsilon$-spaces}
exhaust the class of $n$-dimensional Lorentzian manifolds  admitting a group of isometries of dimension at least $\frac{1}{2} n(n-1)+1$, for almost all values of $n$ [Patrangenaru V.,
\emph{Geom. Dedicata} \textbf{102} (2003), 25--33]. We shall prove that the curvature tensor of these spaces satisfy several interesting algebraic properties. In particular, we will show that Egorov spaces are Ivanov--Petrova manifolds, curvature-Ricci commuting (indeed, semi-symmetric) and $\mathcal P$-spaces, and that $\varepsilon$-spaces are Ivanov--Petrova and curvature-curvature commuting manifolds.}


\Classification{53C50; 53C20}

\section{Introduction}

The study of the relations between algebraic properties of the curvature tensor and the geometry of the underlying manifold,
is a f\/ield of great interest, which has been intensively studied in recent years. In general, one f\/inds dif\/ferent behaviours according to dif\/ferent possibilities for the signature of the metric, and the Riemannian case usually turns out to be much more restrictive than the Lorentzian and the pseudo-Riemannian ones.
In this paper, we shall emphasize some interesting properties, determined by the curvature tensor, for two remarkable classes of $n$-dimensional Lorentzian manifolds, known as {\em Egorov spaces} and {\em $\varepsilon$-spaces}, which naturally occur in the classif\/ication of Lorentzian manifolds admitting a group of isometries of dimension at least $\frac{1}{2} n(n-1)+1$.

These spaces do not have a Riemannian counterpart. In fact, an $n$-dimensional  Riemannian manifold that admits a group of isometries of dimension $\frac 12 n(n-1)+1$, is either of constant sectional curvature or the Riemannian product between an $(n-1)$-dimensional manifold of constant sectional curvature and a line (or circle).  More cases are possible in Lorentzian settings, because of the existence of null submanifolds which are left invariant by group actions.

Let $l_0(n)>l_1(n)>\cdots$ denote the possible dimensions of all groups of isometries of Lorentzian manifolds of dimension $n$. Following \cite{P}, an $n$-dimensional connected Lorentzian manifold $M$ is said to belong to the {\em $j$-stratum} if there is a Lie group $G$, of dimension $l_j(n)$, that acts ef\/fectively on $M$ by isometries. The third stratum is formed by Lorentzian manifolds admitting a group of isometries of dimension at least $\frac 12 n(n-1)+1$.

Lorentzian manifolds belonging to the f\/irst two strata (that is, admitting a group of isometries of dimension at least $\frac 12 n(n-1)+2$) have constant curvature. The complete classif\/ication of Lorentzian manifolds in the third stratum was also given in \cite{P}, in any dimension $n$ greater than $5$ and dif\/ferent from $7$. Besides Lorentzian manifolds of constant curvature $M_1^n(c)$ and manifolds reducible as products $M^{n-1}(c) \times \mathbb R_1$ and $\mathbb R \times M_1^{n-1}(c)$, the remaining examples are:
\begin{itemize}\itemsep=0pt
\item {\em $\varepsilon$-spaces:} Lorentzian manifolds $(\mathbb R^n,g_{\varepsilon})$, where $\varepsilon =\pm 1$ and
\begin{equation}\label{epsi}
g_{\varepsilon} = \sum _{i=1} ^{n-2} (dx_i)^2 - dx_{n-1} dx_n +\varepsilon \left( \sum _{i=1} ^{n-2} x_i ^2\right) (dx_{n-1})^2 .
\end{equation}
\end{itemize}
These spaces are irreducible Lorentzian symmetric spaces which are models for non-symmetric spaces  \cite{CLPTV,P2}.
\begin{itemize}\itemsep=0pt
\item {\em Egorov spaces:} Lorentzian manifolds $(\mathbb R^n,g_{f})$, where $f$ is a positive smooth function of a~real variable and
\begin{equation}\label{ego}
 g_{f} = f(x_n) \sum _{i=1} ^{n-2} (dx_i)^2 +2 dx_{n-1} dx_n.
\end{equation}
\end{itemize}
These manifolds are named after I.P.~Egorov, who f\/irst introduced and studied them in \cite{E}.

The paper is organized in the following way. The basic description
of the curvature of Egorov spaces and $\varepsilon$-spaces, obtained
in \cite{BCD}, will be reported in Section~\ref{section2}. In Sections~\ref{section3}, \ref{section4} and \ref{section5}
we shall respectively deal with properties related to the
skew-symmetric curvature operator, the Jacobi operator and the Szab\'{o}
operator. Among other properties, it is remarkable that both Egorov
spaces and $\varepsilon$-spaces are Ivanov--Petrova
(Theorems~\ref{EgoIP} and~\ref{epsiIP}),  curvature-Ricci and
curvature-curvature commuting as Theorems~\ref{EgoRicciss} and~\ref{epsicc} show. It is worth to point out that Ricci operators are
two-step nilpotent in both cases, as opposed to the case of
simple curvature-Ricci commuting models considered in~\cite{GN2}. In
Section~\ref{section6} we shall prove that Egorov spaces are of recurrent
curvature.  As a consequence, such spaces are also $\mathfrak
P$-spaces. Hence, under this point of view, they are as close as
possible to being locally symmetric, even if most of them are
not even locally homogeneous~\cite{BCD}. We shall conclude with some
general results on Lorentzian manifolds in the third stratum.

 \section[On the curvature of Egorov spaces and $\varepsilon$-spaces]{On the curvature of Egorov spaces and $\boldsymbol{\varepsilon}$-spaces}\label{section2}

In \cite{BCD}, the f\/irst author, W.~Batat and B.~De Leo described
the curvature and the Ricci tensor of Egorov spaces,  investigating
some curvature properties as in particular local symmetry, local
homogeneity and conformal f\/latness.

Let $(\mathbb R^n,g_f)$, $n \geq 3$ denote an Egorov space. As proved in \cite{BCD}, with respect to the basis of coordinate vector f\/ields $\{\partial_i=\frac{\partial}{\partial _{x_i}}\}$ for which \eqref{ego} holds, the possibly non-vanishing covariant derivatives of coordinates vector f\/ields are given by
\begin{equation}\label{nablaei}
\nabla _{\partial_i} \partial _i= -\frac{f'}{2} \partial _{n-1},  \qquad \nabla _{\partial_i} \partial _n =\frac{f'}{2f} \partial _i, \qquad  i=1,\dots,n-2.
\end{equation}
Moreover, the components of the $(0,4)$-curvature tensor,
which is taken with the sign convention
$R(X,Y)=\nabla_{[X,Y]}-[\nabla_X,\nabla_Y]$, and $R(X,Y,Z,T)=\langle
R(X,Y)Z,T\rangle$ are given by
\begin{equation}\label{curvcomp}
R_{inin}= \frac{1}{4f} \left[ (f')^{2}-2ff''\right],  \qquad i=1,\dots,n-2, \qquad
R_{ijkh}=0 \quad {\rm otherwise}
\end{equation}
and the components of the Ricci tensor, ${\rm Ric}(X,Y)=\mbox{\rm
trace}\, \{Z\mapsto R(X,Z)Y\}$, with respect to~$\{
\partial _i \}$ are
 \begin{equation}\label{rhoei}
{\rm Ric}_{nn} = \frac{n-2}{4f^2} [(f')^2-2ff''], \qquad {\rm Ric}_{ij}=0 \quad {\rm otherwise.}
\end{equation}
Consequently, the {\em Ricci operator}  $\widehat{\rm Ric}$, given by
$\langle\widehat{\rm Ric}X,Y\rangle={\rm Ric}(X,Y)$, is described as follows:
\begin{equation}\label{Qei}
{ \widehat{\rm Ric}}=\left(
\begin{array}{cccc}
0 & \cdots & 0 &  0 \\
\vdots & \vdots & \vdots &\vdots \\
0 & \cdots &  0 &  \frac{n-2}{4f^2} [(f')^2-2ff''] \\
0 & \cdots & 0 & 0
\end{array}
\right).
\end{equation}

The following properties were proved in \cite{BCD}:

\begin{theorem}[\cite{BCD}]\label{EgoBCD}
All Egorov spaces $(\mathbb R^n, g_f)$, $n\geq3$,
\begin{itemize}\itemsep=0pt
\item[$(i)$] are geodesically complete;
\item[$(ii)$] admit a parallel null vector field $\partial _{n-1}$;
\item[$(iii)$] have a two-step nilpotent Ricci operator;
\item[$(iv)$] are conformally flat.
\end{itemize}
\end{theorem}

It is also worthwhile to remark that among all Lorentzian
manifolds in the third stratum, Egorov spaces are the only ones
which need not be symmetric (indeed, not even homogeneous). In fact,
as proved in \cite{BCD}, an Egorov space $(\mathbb R^n, g_f)$ is
locally symmetric if and only if its def\/ining function $f$ satisf\/ies
\begin{equation*}
(f')^2-2ff'' =k f^2,
\end{equation*}
where $k$ is a real constant, and is locally homogeneous if and only if either it is locally symmetric, or its def\/ining function  $f$ is a solution of
\begin{equation*}
(f')^2-2ff''=\frac{c_n}{(x_n+d)^2} f^2,
\end{equation*}
for some real constants $c_n \neq 0$ and $d$.

Further observe that, as a consequence of (\ref{nablaei}) the coordinate vector f\/ield $\partial_{n-1}$ is parallel and null, thus showing that the
underlying structure of Egorov spaces is that of a Walker metric~\cite{walker-metrics}.

Next, let $(\mathbb R^n,g_{\varepsilon})$, $n \geq 3$, denote a $\varepsilon$-space as described by~\eqref{epsi}. Then, the Levi-Civita connection of $(\mathbb R^n,g_{\varepsilon})$ is completely determined by
\begin{equation}\label{epsinablaei}
\nabla _{\partial_i} \partial _{n-1}=-2\varepsilon x_i \partial _{n}, \quad i=1,\dots,n-2,  \qquad \nabla _{\partial_{n-1}} \partial _{n-1} =-\varepsilon \sum _{k=1} ^{n-2} x_k \partial _k.
\end{equation}
The curvature components with respect to the coordinate basis $\{\partial_i\}$  are given by
\begin{equation}\label{epsicurvcomp}
R_{in-1in-1}= -\varepsilon,  \quad i=1,\dots,n-2,   \qquad
R_{ijkh}=0 \quad {\rm otherwise}
\end{equation}
(see \cite{BCD}), and the components of the Ricci tensor are
\begin{equation*}
{\rm Ric}_{n-1\, n-1} = -(n-2) \varepsilon, \qquad {\rm  Ric}_{ij}=0 \quad {\rm otherwise.}
\end{equation*}
Moreover, the following results hold:

\begin{theorem}[\cite{BCD, CLPTV}]\label{epsiBCD}
All $\varepsilon$-spaces $(\mathbb R^n, g_{\varepsilon})$, $n\geq3$,
\begin{itemize}\itemsep=0pt
\item[$(i)$] are locally symmetric;
\item[$(ii)$] admit a parallel null vector field $\partial _{n}$;
\item[$(iii)$] have a two-step nilpotent Ricci operator;
\item[$(iv)$] are conformally flat.
\end{itemize}
\end{theorem}

It immediately follows from~(\ref{epsinablaei}) that
$\varepsilon$-spaces are Walker manifolds (indeed pp-waves) with~$\partial_n$ a~parallel null vector f\/ield. We refer to~\cite{Hall} and references therein for more information on the
geometry of pp-waves.

\section{The skew-symmetric curvature operator}\label{section3}

Let $(M,g)$ be a pseudo-Riemannian manifold and $R$ its
curvature tensor. The {\em skew-symmetric curvature operator} of an
oriented non-degenerate $2$-plane $\pi$ of $M$, is def\/ined as
\[
R(\pi)=|g(u,u)g(v,v)-g(u,v)^2|^{-1/2} R(u,v)
\]
and is independent of the oriented pair $\{u,v\}$ of tangent vectors spanning $\pi$. $(M,g)$ is said to be an {\em Ivanov--Petrova manifold} (shortly, an {\em IP manifold}) if the eigenvalues of $R(\pi)$ are constant on the Grassmannian $G^+(2,n)$ of all oriented $2$-planes. The eigenvalues may change at dif\/ferent points of $(M,g)$. These manifolds are named after S.~Ivanov and I.~Petrova, who introduced them and f\/irst undertook their study~\cite{IP}.

Riemannian IP manifolds have been completely classif\/ied in all dimensions $n\geq 4$ (see \mbox{\cite{GLS,G,N}}). Examples of three-dimensional Riemannian IP manifolds which are neither conformally f\/lat nor curvature homogeneous may be found in~\cite{C} and~\cite{N}.

Recently, IP manifolds have been extended and largely investigated
in pseudo-Riemannian settings  (see
\cite{GRH,G1,GZ} and references therein).

Consider now an Egorov space $(\mathbb R^n, g_f)$.  Let $\pi$ denote
any oriented non-degenerate $2$-plane and $\{u,v\}$ an orthonormal
basis of $\pi$ and put $\Xi=|g(u,u)g(v,v)-g(u,v)^2|^{1/2}$. The
metric components with respect to $\{\partial_i \}$ can be deduced
at once from~\eqref{ego}, and the curvature components are given by
\eqref{curvcomp}. Then, a direct calculation shows that the
skew-symmetric curvature operator  $R(\pi)=\frac{1}{\Xi}
R(u,v)$ takes the form
\begin{gather}\label{skewEgo}
R(\pi)\! = \! \frac{(f')^{2}-2ff''}{4\Xi f^2} \!\left( \!\!
\begin{array}{@{}c@{\,\,\,}c@{\,\,\,}c@{\,\,\,}c@{\,\,\,}c@{}}
0&0&\hdots&0&u_{n}v_1-u_1v_{n}\\
\vdots&\vdots&\vdots&\vdots&\vdots\\
0&0&\hdots&0&u_nv_{n-2}-u_{n-2}v_n\\
f(u_1v_n-u_nv_1)&\hdots&f(u_{n-2}v_n-u_nv_{n-2})&0&0\\
0&0&\hdots&0&0
\end{array}
\!\!\right)\!,\!\!
\end{gather}
where $u=\sum u_i \partial _i$ and $v=\sum v_i \partial _i$.
It is easily seen, by \eqref{skewEgo}, that the eigenvalues of the
skew-symmetric curvature operator vanish identically. This can be
proved either by calculating the eigenvalues, or as a consequence of
the fact that $R(\pi)$ is three-step nilpotent. Since all the
eigenvalues of the skew-symmetric curvature operator are constant,
we proved the following

\begin{theorem}\label{EgoIP}
All Egorov spaces are IP manifolds.
\end{theorem}

We now turn our attention to $\varepsilon$-spaces. Using
\eqref{epsi} and \eqref{epsicurvcomp}, it is easy to prove that the
skew-symmetric curvature operator associated to any oriented
non-degenerate $2$-plane $\pi$ is given by
\begin{gather}\label{skewepsi}
R(\pi) =  \frac{\varepsilon}{\Xi}\!  \left(
\begin{array}{@{}c@{\,\,\,}c@{\,\,\,}c@{\,\,\,}c@{\,\,\,}c@{}}
0 & \hdots & 0 & u_1v_{n-1}-v_1u_{n-1}& 0  \\
0 & \hdots & 0   & u_2v_{n-1}-v_2u_{n-1}& 0\\
\vdots & \vdots  & \vdots & \vdots & \vdots\\
0 & \hdots & 0 & u_{n-2}v_{n-1}-v_{n-2}u_{n-1}& 0  \\
0 & \hdots & 0 & 0 & 0  \\
u_1v_{n-1}-v_1u_{n-1} & \hdots &u_{n-2}v_{n-1}-v_{n-2}u_{n-1} & 0& 0
\end{array}
 \right)\!,\!\!
\end{gather}
where $\{u,v\}$, $u=\sum u_i \partial _i$, $\, v=\sum v_i \partial
_i$, is a positively oriented orthonormal basis of $\pi$ and
$\Xi=|g(u,u)g(v,v)-g(u,v)^2|^{1/2}$.
 Now, \eqref{skewepsi} implies that the eigenvalues of the
skew-symmetric curvature operator vanish identically (in particular,
$R(\pi)$ is again three-step nilpotent). Hence, we have the
following

\begin{theorem}\label{epsiIP}
All $\varepsilon$-spaces are IP manifolds.
\end{theorem}

 \section{The Jacobi and Szab\'{o} operators}\label{section4}

Let $(M,g)$ be a pseudo-Riemannian manifold, $p$ a point of $M$ and $u \in T_p M$ a unit vector. Denote by $\gamma : r \mapsto \exp (ru)$ the geodesic with arc length $r$ tangent to $u$ at $p$. The {\em Jacobi operator} along $\gamma$ is def\/ined by  $R_{\gamma} = R(\dot{\gamma} , \cdot ) \dot{\gamma} $ (here, we are adopting the sign convention $R(u,v)=\nabla_{[u,v]}-[\nabla_u,\nabla_v]$ for the curvature tensor). In particular, $\mathcal J _u =R_{\gamma}(0)$ def\/ines the Jacobi operator at $u$. The  {\em Szab\'{o} operator} at $u$ is def\/ined by $\mathcal S _u := \mathcal J ' _u := D_{\dot{\gamma}} R_{\gamma} (0)$.
The role of the Jacobi operator in determining the geometry of $(M,g)$ has been extensively studied.
In particular, two natural generalizations of locally symmetric spaces were introduced in \cite{BV} in terms of the Jacobi operator. A manifold $(M,g)$ is said to be
\begin{itemize}\itemsep=0pt
\item a {\em $\mathfrak C$-space} if the eigenvalues of
$R_{\gamma}$ are constant along $\gamma$, for each geodesic $\gamma$ of $(M,g)$;
\item a {\em $\mathfrak P$-space} if the Jacobi operator along any geodesic can be diagonalized by a parallel orthonormal frame along the geodesic.
\end{itemize}
As proved in \cite{BV, CL-GR-VA-VL}, a Riemannian or Lorentzian manifold $(M,g)$ is locally symmetric if and only if it is both a $\mathfrak C$-space and a $\mathfrak P$-space. Moreover, at least in the analytic case,
$(M,g)$ is a $\mathfrak P$-space if and only if
\begin{equation}\label{JJ'}
[\mathcal S _u,\mathcal J _u]=0 \qquad \text{for all} \; u,
\end{equation}
that is, when the Jacobi and Szab\'{o} operators of all tangent vectors $u$ commute.
However, these results are no longer true in higher signatures, where examples exist which are $\mathfrak C$- and $\mathfrak P$-spaces simultaneously but
not locally symmetric.
Now, being symmetric, a $\varepsilon$-space is obviously both a~$\mathfrak C$-space and a~$\mathfrak P$-space.  As concerns Egorov spaces, we calculated the Jacobi operator starting from~\eqref{curvcomp}, obtaining
\begin{equation}\label{JacEgo}
\mathcal J _u=  (n-2)\frac{(f')^{2}-2ff''}{4f^2} \left(
\begin{array}{cccccc}
u_n^2&0&\hdots&\hdots&0&-u_1u_n\\
0&u_n^2&0&\hdots&0&-u_2u_n\\
\vdots&\vdots&\vdots&\vdots&0&\vdots\\
0&0&\hdots&u_n^2&0&-u_{n-2}u_n\vspace{1mm}\\
-fu_1u_n&-fu_2u_n&\hdots&-fu_{n-2}u_n&0&f\sum\limits_{i\leq n-2} u_i^2\vspace{1mm}\\
0&0&0&0&0&0
\end{array}
\right),
\end{equation}
where $u=\sum u_i \partial _i$, while for the Szab\'{o} operator,
from \eqref{nablaei} and \eqref{JacEgo}, we get
\begin{equation}\label{SzEgo}
\mathcal S _X=  \frac{1}{u_n}\left(\frac{(f')^{2}-2ff''}{4f^2}\right)'\mathcal J _u.
\end{equation}
From \eqref{JacEgo} and \eqref{SzEgo} it follows at once that
\eqref{JJ'} holds for any tangent vector $x$. Thus, we proved the following

\begin{theorem}\label{EgoP}
All Egorov spaces are $\mathfrak P$-spaces.
\end{theorem}

We f\/inish this section with the following remark on
$\mathfrak{C}$-spaces. Since the eigenvalues of the Jacobi operators
$\mathcal J_u$ are completely determined by the traces of the powers
$\mathcal J^k_u$ for all $k$, a~Lorentzian manifold is a
$\mathfrak{C}$-space if and only if $\nabla_{\gamma'}\,\mbox{\rm
trace}\, \mathcal{J}_{\gamma}^k=0$ for all $k$ and all geodesics
$\gamma$.

The f\/irst of those conditions above ($\nabla_{\gamma'}\,\mbox{\rm
trace}\, \mathcal{J}_{\gamma}=0$) shows that the Ricci tensor is
cyclic parallel (i.e., $(\nabla_X{\rm Ric})(X,X)=0$ for all vector f\/ields
$X$). Now, a straightforward calculation using the expressions in
Section~\ref{section2} shows that an Egorov space has cyclic parallel Ricci
tensor if and only if it is locally symmetric (i.e.,
$(f')^2-2ff''=4f^2 C$ for some constant $C$).

 \section{Commuting curvature operators and semi-symmetry}\label{section5}

Commuting properties of the curvature, Ricci and Jacobi operators have been intensively studied in the last years. For several of these properties, a complete description of the corresponding manifolds has not been obtained yet. A recent survey can be found in~\cite{MFGGV}, and we can refer to~\cite{GRH1} for the sistematic study of the three-dimensional Lorentzian case.

Let $(M,g)$ be a pseudo-Riemannian manifold and denote by
$\widehat{\rm Ric}$ and $\mathcal J$ the Ricci operator and the Jacobi
operator of $(M,g)$, respectivey. The manifold is said to be
\begin{itemize}\itemsep=0pt
\item {\em Jacobi--Ricci commuting} if $\mathcal J _u \cdot {
\widehat{\rm Ric}} = { \widehat{\rm Ric}} \cdot \mathcal J _u$ for any tangent
vector $u$;
\item {\em curvature-Ricci commuting} if $R(u,v) \cdot {
\widehat{\rm Ric}} = { \widehat{\rm Ric}} \cdot R (u,v)$ for all tangent vectors
$u$, $v$;
\item {\em curvature-curvature commuting} if $R(u,v) \cdot R(z,w) = R(z,w) \cdot R (u,v)$ for all tangent vectors $u$, $v$, $z$, $w$.
\end{itemize}
The class of Jacobi--Ricci commuting manifolds coincides with the one of curvature-Ricci commuting manifolds~\cite{GN}, which are also known in literature as {\em Ricci semi-symmetric spaces}.
(We may refer to~\cite{HV} for geometric interpretations and further references.)
Other commuting curvature conditions, related to the {\em Weyl conformal curvature tensor}~$W$, are also well known~\cite{MFGGV}. However, since both Egorov spaces and $\varepsilon$-spaces are conformally f\/lat, their conformal curvature tensor identically vanishes and so, the commuting curvature conditions involving $W$ are trivially satisf\/ied.

For Egorov spaces $(\mathbb R^n,g_f)$, using \eqref{Qei} and~\eqref{skewEgo}, a straightforward calculation leads to a~proof of the following

\begin{theorem}\label{EgoRicciss}
All Egorov spaces are curvature-Ricci commuting and curvature-curvature commuting.
\end{theorem}

We brief\/ly recall that a pseudo-Riemannian manifold $(M,g)$ is said to be {\em semi-symmetric} if its curvature tensor $R$ satisf\/ies $R (u, v) \cdot R = 0,$
for all tangent vectors~$u$,~$v$. Semi-symmetric spaces are a well known and natural generalization of locally symmetric spaces. A semi-symmetric space is Ricci semi-symmetric. The converse does not hold in general, but it is true for conformally f\/lat manifolds, since their curvature is completely determined by the Ricci tensor.  Hence, from Theorem \ref{EgoRicciss} and $(iv)$ of Theorem \ref{EgoBCD} we have at once the following

\begin{theorem}\label{Egoss}
All Egorov spaces are semi-symmetric.
\end{theorem}
Since $\varepsilon$-spaces are symmetric, they are semi-symmetric and thus Ricci semi-symmetric (i.e., curvature-Ricci commuting). Moreover, using~\eqref{skewepsi}, a direct calculation proves the following

\begin{theorem}\label{epsicc}
All $\varepsilon$-spaces are curvature-curvature commuting.
\end{theorem}

 \section{Recurrent curvature}\label{section6}

A pseudo-Riemannian manifold $(M,g)$ is said to be of {\em recurrent curvature} if the covariant derivative of its curvature tensor is linearly dependent of the curvature tensor itself. As Egorov spaces $(\mathbb R^n,g_f)$ are concerned, with respect to coordinate vector f\/ields $\{ \partial_i \}$, by \eqref{rhoei} we have
\begin{equation}\label{ro}
\varrho= (n-2)\frac{(f')^2-2ff''}{4f^2}
\left(\begin{array}{cccc}
0 & \dots  & 0& 0 \\
\vdots &  \vdots & \vdots \\
0 & \dots  & 0 & 0\\
0 & \dots & 0 & 1
\end{array}\right)
\end{equation}
and so, \eqref{nablaei} and \eqref{ro} easily give
\begin{equation*}
\nabla \varrho= (n-2)\left(\frac{(f')^2-2ff''}{4f^2}\right)'
\left(\begin{array}{cccc}
0 & \dots & 0& 0 \\
\vdots &  \vdots & \vdots \\
0 & \dots  & 0 & 0\\
0 & \dots & 0 & 1
\end{array}\right).
\end{equation*}
Thus, the Ricci tensor of $(\mathbb R^n,g_f)$ is recurrent. Since
$(\mathbb R^n,g_f)$ is conformally f\/lat by $(iv)$ of Theorem \ref{EgoBCD}, we have at once the following

\begin{theorem}\label{recEgo}
All Egorov spaces are of recurrent curvature.
\end{theorem}

\section{General conclusions}\label{section7}

As we already mentioned in the Introduction, together with Lorentzian spaces of constant curvature and Lorentzian products of a space of constant curvature and a real line (or circle), Egorov spaces and $\varepsilon$-spaces complete the classif\/ication of Lorentzian manifolds in the third stratum, in all dimensions $n>5$, $n\neq 7$. The socalled $\mathfrak h$-triple method, used in \cite{P} to obtain the classif\/ication, does not apply for $n=7$, even if Egorov and $\varepsilon$-spaces are def\/ined in all dimensions greater than two. The classif\/ication above, together with the results of the previous sections, lead to the following.
\begin{quote}
\emph{For any dimension $n>5$, $n\neq 7$, Lorentzian manifolds $(M^n,g)$ in the third stratum are locally conformally flat and moreover,
they are curvature-Ricci commuting $($and thus semi-symmetric$)$ and of recurrent curvature $($and thus $\mathfrak P$-spaces$)$.}
\end{quote}
Some of the properties listed above are trivial for all Lorentzian manifolds in the third stratum except for Egorov spaces. However, they are interesting since they f\/ix a common base of properties for manifolds belonging to the third stratum. As shown in \cite{BCD}, other properties, as symmetry and homogeneity, are not the common inheritance for these manifolds.
We refer to Theorem~5.3 of~\cite{BCD} for other properties shared by all Lorentzian manifolds in the third stratum.

As concerns Lorentzian manifolds in the third stratum without Riemannian counterpart, we have in addition the following
\begin{quote}
\emph{Egorov-spaces and $\varepsilon$-spaces are Ivanov--Petrova and curvature-curvature commuting Walker manifolds}.
\end{quote}

\subsection*{Acknowledgements}

First author supported by funds of MIUR (PRIN
2007) and the University of Salento. Second author supported by
projects MTM2009-07756 and INCITE09 207 151 PR  (Spain).
Finally the authors would like to express their thanks to the Referees of this paper for pointing out some
mistakes in the original manuscript.

\pdfbookmark[1]{References}{ref}
\LastPageEnding


\begin{thebibliography}{99}

\footnotesize\itemsep=0pt

\bibitem{BCD}
Batat W., Calvaruso G., De Leo B.,
Curvature properties of Lorentzian manifolds with large isometry groups,
\href{http://dx.doi.org/10.1007/s11040-009-9060-4}{\emph{Math. Phys. Anal. Geom.}} \textbf{12} (2009), 201--217.

\bibitem{BV}
Berndt J., Vanhecke L.,
Two natural generalizations of locally symmetric spaces,
\href{http://dx.doi.org/10.1016/0926-2245(92)90009-C}{\emph{Differential Geom. Appl.}} \textbf{2} (1992), 57--80.

\bibitem{MFGGV}
Brozos-V\'{a}zquez~M., Fiedler~B., Garc\'{i}a-R\'{i}o~E., Gilkey~P., Nik\v{c}evi\'{c}~S., Stanilov~G., Tsankov~Y., V\'{a}zquez-Lorenzo~R., Videv~V.,
Stanilov--Tsankov--Videv theory,
\href{http://dx.doi.org/10.3842/SIGMA.2007.095}{\emph{SIGMA}} \textbf{3} (2007), 095, 13~pages,
\href{http://arxiv.org/abs/0708.0957}{arXiv:0708.0957}.


\bibitem{walker-metrics}
Brozos-V\'{a}zquez~M., Garc\'{i}a-R\'{i}o~E., Gilkey P., Nik\v{c}evi\'{c}~S., V\'{a}zquez-Lorenzo~R.,
The geometry of Walker manifolds,
{\it Synthesis Lectures on Mathematics and Statistics}, Vol.~5, Morgan \& Claypool Publ., 2009.



\bibitem{C}
Calvaruso G.,
Three-dimensional Ivanov--Petrova manifolds,
\href{http://dx.doi.org/10.1063/1.3152607}{\emph{J. Math. Phys.}} \textbf{50} (2009), 063509, 12~pages.

\bibitem{CL-GR-VA-VL}
Calvi\~no-Louzao E., Garc\'{\i}a-R\'{\i}o E., V\'{a}zquez-Abal M.E., V\'{a}zquez-Lorenzo~R.,
Curvature operators and ge\-ne\-ralizations of  symmetric spaces in Lorentzian geometry, Preprint, 2009.

\bibitem{CLPTV}
Cahen M., Leroy J., Parker M., Tricerri F., Vanhecke L.,
Lorentz manifolds modelled on a Lorentz symmetric space,
\href{http://dx.doi.org/10.1016/0393-0440(90)90007-P}{\emph{J. Geom. Phys.}} \textbf{7} (1990), 571--581.

\bibitem{E}
Egorov I.P.,
Riemannian spaces of the f\/irst three lacunary types in the geometric sense,
\emph{Dokl. Akad. Nauk. SSSR} \textbf{150} (1963), 730--732.

\bibitem{GRH}
Garc\'{i}a-R\'{i}o E., Haji-Badali A., V\'{a}zquez-Lorenzo R.,
Lorentzian three-manifolds with special curvature ope\-ra\-tors,
\href{http://dx.doi.org/10.1088/0264-9381/25/1/015003}{\emph{Classical Quantum Gravity}} \textbf{25} (2008), 015003, 13~pages.

\bibitem{GRH1}
Garc\'{i}a-R\'{i}o E., Haji-Badali A., V\'{a}zquez-Abal M.E., V\'{a}zquez-Lorenzo~R.,
Lorentzian $3$-manifolds with commuting curvature operators,
\href{http://dx.doi.org/10.1142/S0219887808002941}{\emph{Int. J. Geom. Methods Mod. Phys.}} \textbf{5} (2008), 557--572.

\bibitem{G}
Gilkey P.,
Riemannian manifolds whose skew-symmetric curvature operator has constant eigenvalues.~II, in
Dif\/ferential Geometry and Applications (Brno, 1998), Masaryk Univ., Brno, 1999, 73--87.

\bibitem{G1}
Gilkey P.,
Geometric properties of natural operators def\/ined by the Riemann curvature tensor, World Scien\-tif\/ic Publishing Co. Inc., River Edge, NJ, 2001.

\bibitem{GLS}
Gilkey P., Leahy J.V., Sadofsky H.,
Riemannian manifolds whose skew-symmetric curvature operator has constant eigenvalues,
\href{http://dx.doi.org/10.1512/iumj.1999.48.1699}{\emph{Indiana Univ. Math. J.}} \textbf{48} (1999), 615--634.

\bibitem{GN}
Gilkey P., Nik\v{c}evi\'{c} S.,
Pseudo-Riemannian Jacobi--Videv manifolds,
\href{http://dx.doi.org/10.1142/S0219887807002272}{\emph{Int. J. Geom. Methods Mod. Phys.}} \textbf{4} (2007), 727--738,
\href{http://arxiv.org/abs/0708.1096}{arXiv:0708.1096}.

\bibitem{GN2}
Gilkey P., Nik\v{c}evi\'{c} S.,
The classif\/ication of simple Jacobi--Ricci commuting algebraic curvature tensors,
\emph{Note Mat.} \textbf{28} (2008), suppl.~1, 341--348,
\href{http://arxiv.org/abs/0710.2080}{arXiv:0710.2080}.

\bibitem{GZ}
Gilkey P., Zhang T.,
Algebraic curvature tensors for indef\/inite metrics whose skew-symmetric curvature operator has constant Jordan normal form,
\emph{Houston J. Math.} \textbf{28} (2002), 311--328,
\href{http://arxiv.org/abs/math.DG/0205079}{math.DG/0205079}.

\bibitem{HV}
Haesen S., Verstraelen L.,
Natural intrinsic geometrical symmetries,
\href{http://dx.doi.org/10.3842/SIGMA.2009.086}{{\em SIGMA}} {\bf 5} (2009),  086, 14 pages,
\href{http://arxiv.org/abs/0909.0478}{arXiv:0909.0478}.


\bibitem{Hall}
 Hall G.S.,
 Symmetries and curvature structure in general relativity,
 {\it World Scientific Lecture Notes in Physics}, Vol.~46, World Scientif\/ic Publishing Co., Inc., River Edge, NJ, 2004.

\bibitem{IP}
Ivanov S., Petrova I.,
Riemannian manifolds in which certain curvature operator has constant eigenvalues along each circle,
\href{http://dx.doi.org/10.1023/A:1006548328030}{\emph{Ann. Global Anal. Geom.}} \textbf{15} (1997), 157--171.

\bibitem{N}
Nikolayevsky Y.,
Riemannian manifolds whose curvature operator $R(X,Y)$ has constant eigenvalues,
\href{http://dx.doi.org/10.1017/S0004972700034523}{\emph{Bull. Austral. Math. Soc.}} \textbf{70} (2004), 301--319, \href{http://arxiv.org/abs/math.DG/0311429}{math.DG/0311429}.

\bibitem{P}
Patrangenaru V.,
Lorentz manifolds with the three largest degrees of symmetry,
\href{http://dx.doi.org/10.1023/B:GEOM.0000006588.95481.1c}{\emph{Geom. Dedicata}} \textbf{102} (2003), 25--33.

\bibitem{P2}
Patrangenaru V.,
Locally homogeneous pseudo-Riemannian manifolds,
\href{http://dx.doi.org/10.1016/0393-0440(94)00040-B}{\emph{J. Geom. Phys.}} \textbf{17} (1995), 59--72.

\end{thebibliography}
\end{document}